\documentclass[12pt]{article}
\usepackage{amsthm,amssymb,latexsym}
\title{Construction of doubly-periodic instantons}
\author{Marcos Jardim \\ Yale University \\ Department of Mathematics \\ 10 Hillhouse Avenue \\
New Haven, CT 06520-8283 USA}

\newcommand{\pf}{{\em Proof: }} \newcommand{\pfend}{\hfill $\Box$ \linebreak}
\newcommand{\seta}{\rightarrow} \newcommand{\torus}{T\times\cpx}
\newcommand{\tproj}{T\times\proj} 
 \newcommand{\dual}{\hat{T}}
 \newcommand{\del}{\overline{\partial}}

\newcommand{\ksi}{\xi} 
\newcommand{\real}{\mathbb{R}} \newcommand{\cpx}{\mathbb{C}} \newcommand{\zed}{\mathbb{Z}}
\newcommand{\proj}{\mathbb{P}^1} \newcommand{\vv}{{\cal V}}
\newcommand{\ee}{{\cal E}}  \newcommand{\pp}{{\cal P}}
  
  \newcommand{\jj}{{\cal J}}
\newcommand{\ind}{{\rm index}}  
 \newcommand{\as}{\pm\ksi_0}
\newcommand{\oo}{{\cal O}}  \newcommand{\rr}{{\cal R}}

\newtheorem{thm}{Theorem}
\newtheorem{lem}[thm]{Lemma}
\newtheorem{prop}[thm]{Proposition}
\newtheorem{conj}[thm]{Conjecture}

\begin{document}
\maketitle

\begin{abstract}
We construct finite-energy instanton connections over $\real^4$
which are periodic in two directions via an analogue of the Nahm
transform for certain singular solutions of Hitchin's equations
defined over a 2-torus.
\end{abstract}

\baselineskip18pt \newpage

\section{Introduction} \label{intro}

Since the appearance of the Yang-Mills equation on the mathematical
scene in the late 70's, its anti-self-dual (ASD) solutions have
been intensively studied. The first major result in the field was the
ADHM construction of instantons on $\real^4$ \cite{ADHM}. Soon
after that, W. Nahm adapted the ADHM construction to obtain the
{\em time-invariant} ASD solutions of the Yang-Mills equations, the
so-called monopoles \cite{N}. It turns out that these constructions
are two examples of a much more general framework.

The {\em Nahm transform} can be defined in general for anti-self-dual
connections on $\real^4$, which are invariant under some subgroup of 
translations $\Lambda\subset\real^4$ (see \cite{Nak2}). In these generalised
situations, the Nahm transform gives rise to {\em dual instantons} on 
$(\real^4)^*$, which are invariant under:
$$ \Lambda^*=\{\alpha\in(\real^4)^*\ |\ \alpha(\lambda)\in\zed\ \forall\lambda\in\Lambda\} $$
There are plenty of examples of such constructions available in
the literature, namely:

\begin{itemize}
\item The trivial case $\Lambda=\{0\}$ is closely related to the
celebrated ADHM construction of instantons, as described by Donaldson
\& Kronheimer \cite{DK}; in this case, $\Lambda^*=(\real^4)^*$ and
an instanton on $\real^4$ corresponds to some algebraic data.

\item If $\Lambda=\zed^4$, this is the Nahm transform of Braam \& van
Baal \cite{BVB} and Donaldson \& Kronheimer \cite{DK}, defining a
hyperk\"ahler isometry of the moduli space of instantons over two dual
4-tori.

\item $\Lambda=\real$ gives rise to monopoles, extensively studied by
Hitchin \cite{H3}, Donaldson \cite{D}, Hurtubise \& Murray \cite{HM} 
and Nakajima \cite{Nak2}, among several others; here, $\Lambda^*=\real^3$, 
and the transformed object is, for SU(2) monoples, an analytic solution of
certain matrix-valued ODE's (the so-called Nahm's equations), defined
over the open interval $(0,2)$ and with simple poles at the end-points.

\item $\Lambda=\zed$ correspond to the so-called calorons, studied
by Nahm \cite{N}, Garland \& Murray \cite{GM} and others; the
transformed object is the solution of certain nonlinear Nahm-type
equations on a circle.

\end{itemize}

The purpose of this paper fits well into this larger mathematical
programme. Our goal is to construct finite-energy instantons over
$\real^4$ provided with the Euclidean metric, which are periodic in
two directions ($\Lambda^*=\zed^2$), so-called {\em doubly-periodic
instantons}, from solutions of Hitchin's equations \cite{H} defined
on a 2-torus, i.e. instantons over $\real^4$ which are invariant
under $\Lambda=\zed^2\times\real^2$. The latter object is now very
well studied, and their existence is determined by certain
holomorphic data.

One might also ask if all doubly-periodic instantons can be produced
in this way. In the sequel \cite{J2} of this paper, we will show that
the construction here presented is invertible by describing the Nahm
transform for instantons over $T^2\times\real^2$, which produce
singular solutions of Hitchin's equations.

Indeed, Hitchin's equations admit very few smooth solutions over
elliptic curves (see \cite{H}). Therefore, by analogy with
Hitchin's construction of monopoles \cite{H3}, we will consider a
certain class of singular solutions, for which existence is
guaranteed \cite{Ko,S}. The singularity data is converted into the
asymptotic behaviour of the Nahm transformed doubly-periodic
instanton; such a picture is again familiar from the construction
of monopoles.

A string-theoretical version of the Nahm transform here presented
was given by Kapustin \& Sethi \cite{KS}. In fact, the other examples
of Nahm transforms mentioned above also have string-theoretical
interpretations. The ADHM construction and the Fourier transform of
instantons over 4-tori were discussed in these terms by Witten \cite{W},
while Kapustin \& Sethi \cite{KS} also treated the case of calorons.

Let us now outline the contents of this paper. Section \ref{higgs}
is dedicated to a brief review of Hitchin's self-duality equations,
and the precise description of the particular type of solutions we
will be interested in. The main topic of the paper is contained in
sections \ref{inv} and \ref{invholo}, when we will show how to
construct doubly-periodic instantons and explore some of the
properties of the instantons obtained. We conclude with a few
remarks and raising some questions for future investigation.

\paragraph{Acknowledgements.}
This work is part of my Ph.D. project \cite{J}, which was funded by
CNPq, Brazil. I am grateful to my supervisors, Simon Donaldson and
Nigel Hitchin, for their constant support and guidance. I also
thank Olivier Biquard, Alexei Kovalev and Brian Steer for
invaluable help in the later stages of the project.


\section{Singular Higgs pairs} \label{higgs}

In \cite{H} Hitchin studied the dimensional reduction of the usual
Yang-Mills anti-self-dual equations from four to two dimensions. More
precisely, let $V\seta\real^4$ be a rank $k$ vector bundle with a
connection $\tilde{B}$ which does not depend on two coordinates. Pick
up a global trivialisation of $V$ and write down $\tilde{B}$ as
a 1-form:
$$ \tilde{B}=B_1(x,y)dx+B_2(x,y)dy+\phi_1(x,y)dz+\phi_2(x,y)dw $$
Hitchin then defined a {\em Higgs field} $\Phi=(\phi_1-i\phi_2)d\xi$,
where $d\xi=dx+idy$. So $\Phi$ is a section of $\Lambda^{1,0}{\rm End}V$,
where $V$ is now seen as a bundle over $\real^2$ with a connection
$B=B_1dx+B_2dy$.

The ASD equations for $\tilde{B}$ over $\real^4$ can then be
rewitten as a pair of equations on $(B,\Phi)$ over $\real^2$:
\begin{equation} \label{hiteqs} \left\{ \begin{array}{l}
F_B+[\Phi,\Phi^*]=0 \\
\del_B\Phi=0
\end{array} \right. \end{equation}
These equations are also conformally invariant, so they make sense
over any Riemann surface. Solutions $(B,\Phi)$ are often called {\em
Higgs pairs}.

As we mentioned in the introduction, we are interested in singular
Higgs pairs over a 2-torus $\dual$ defined on an $U(k)$-bundle
$V\seta\dual$. Since we want to think of $\dual$ as a quotient of
$\real^4$ by $\Lambda=\zed^2\times\real^2$, the natural choice of
metric for $\dual$ is the flat, Euclidean metric. Let us also fix a
complex structure on $\dual$ coming from a choice of complex
structure on $\real^4$.

Singular Higgs bundles were widely studied by many authors
(\cite{S}, \cite{K} and \cite{Ko} among others) and are closely
related to the so-called {\em parabolic Higgs bundles}. Adopting
this point of view, we will consider a holomorphic vector bundle
$\vv\seta\dual$ of degree $-2$ with the following quasi-parabolic
structure over two points $\as\in\dual$ (regarding now $\dual$ as
an elliptic curve):
\begin{eqnarray*}
\begin{array}{rcl}
\vv_{\as} = F_1\vv_{\as} \supset &
\underbrace{F_2\vv_{\as}} & \supset F_3\vv_{\as} =\{0\} \\
& {\rm dim}=1 &
\end{array} &
& {\rm order}(\ksi_0) \neq 2
\\ & & \\
\begin{array}{rccl}
\vv_{\xi_0} = F_1\vv_{\xi_0} \supset & \underbrace{F_2\vv_{\xi_0}} & 
\supset \ \underbrace{F_3\vv_{\xi_0}} & \supset F_4\vv_{\xi_0} =\{0\} \\
& {\rm dim}=2 & \ \ \ {\rm dim}=1 &
\end{array} &
& {\rm order}(\ksi_0) = 2
\end{eqnarray*}
To complete the parabolic structure we need to assign {\em weights}
$\alpha_1(\as)$ to $F_1\vv_{\as}$ and $\alpha_2(\as)$ to $F_2\vv_{\as}$
if $\xi_0\neq-\xi_0$ or $\alpha_1(\xi_0)$ to $F_1\vv_{\xi_0}$, $\alpha_2(\xi_0)$ 
to $F_2\vv_{\xi_0}$ and $\alpha_3(\xi_0)$ to $F_3\vv_{\xi_0}$ if $\xi_0=-\xi_0$.
We assume that $\alpha_1=0$ in both cases; if $\xi_0$ is not of order two,
we fix that $\alpha_2(\xi_0)=1+\alpha$ and $\alpha_2(-\xi_0)=1-\alpha$; if 
$\xi_0$ has order two, we fix that $\alpha_2(\xi_0)=1-\alpha$ and 
$\alpha_3(\xi_0)=1+\alpha$ for some $0\leq\alpha<\frac{1}{2}$.
Note in particular that $\vv$ with this parabolic structure has zero
parabolic degree.

From the point of view of the Higgs pair $(B,\Phi)$, this means
that the bundle $V$ is defined away from $\as$, and satisfies,
holomorphically:
$$ \vv|_{\dual\setminus\{\as\}} \simeq (V,\del_B) $$
The Higgs field $\Phi$ has simple poles at the parabolic points
$\as\in\dual$ such that the residues $\phi_0(\as)$ of $\Phi$ are
$k\times k$ matrices  of rank 1. If $\xi_0$ is one of the four
elements of order 2 in $\dual$, then the residue $\phi_0(\xi_0)$ is
assumed to be a $k\times k$ matrix of rank 2. 

Moreover, the harmonic metric $h$ associated with the Higgs pair 
$(B,\Phi)$ is assumed to be compatible with the parabolic structure. 
This means that, in a holomorphic trivialisation of $V$ over a 
sufficiently small neighbourhood around $\as$, $h$ is non-degenerate 
along the kernel of the residues of $\Phi$, and $h\sim O(r^{1\pm\alpha})$ 
along the image of the residues of $\Phi$.

Such metric is clearly not a hermitian metric on the extended 
bundle $\vv$ (since it degenerates at $\as$). Let $h'$ be a
hermitian metric on $\vv$ bounding above the harmonic metric
on $V$.

If $(\vv,\Phi)$ is $\alpha$-stable in the sense of parabolic Higgs
bundles, then the existence of a meromorphic Higgs pair as above is
guaranteed \cite{S} for any rank $k$ and any choice of $\as$. 

Moreover, one usually fixes the eigenvalues of the residues of 
$\Phi$ as well. In our situation, this amounts to choosing only 
one complex number that we denote by $\epsilon$. We assume that 
$\epsilon\neq0$, i.e. the residues of $\Phi$ are semi-simple.

However, in this paper, these parameters (the weights $\alpha_i$ 
and the eigenvalue of the residues $\epsilon$) will be allowed to 
vary; see \cite{BJ} for a complete discussion. It is reassuring to 
know that if two sets of parameters $(\alpha,\epsilon)$ and 
$(\alpha',\epsilon')$ are chosen in generic position, then 
$\alpha$-stability and $\alpha'$-stability are in fact equivalent
conditions \cite{Nak}.

In particular, the case $k=1$ is very simple: once the parameters
parameters $(\alpha,\epsilon)$ are fixed and for any choice of $\as$, 
the moduli space of meromorphic Higgs pairs is just the cotangent 
bundle of $T$, that is a copy of $\torus$.

We will study solutions of (\ref{hiteqs}) over $\dual$
with the singularities $\as$ removed. Due to the non-compactness
of $\dual\setminus\{\as\}$, the choice of metric on the base space
is a delicate issue. From the point of view of the Nahm transform,
it is important to consider the Euclidean, incomplete metric on
the punctured torus, as it is well-known from the examples mentioned
above. However, such a choice of metric is not a good one from the
analytical point of view. For instance, one cannot expect, on general
grounds, to have a finite dimensional moduli space of Higgs pairs.

Fortunately, as we mentioned before, Hitchin's equations are
conformally invariant, so that we are allowed to make conformal
changes in the Euclidean metric localised around the punctures to
obtain a complete metric on \linebreak $\dual\setminus\{\as\}$.
Thus, our strategy is to obtain results concerning the Euclidean
metric from known statements about complete metrics.

In \cite{B3}, Biquard considered the so-called {\em Poincar\'e metric},
which is defined as follows. We perform a conformal change on the
incomplete metric over the punctured torus localised on small punctured
neighbourhoods $D_0$ of $\as$, so that if $\ksi=(r,\theta)$ is a local
coordinate on $D_0$, we have the metric:
\begin{equation} \label{pmetric}
ds^2_P \ = \ \frac{d\ksi d\overline{\ksi}}{|\ksi|^2\log^2|\ksi|^2} \
= \ \frac{dr^2}{r^2\log^2r}+\frac{d\theta^2}{4\log^2r}
\end{equation}
We denote the complete metric so obtained by $g_P$. The Euclidean
metric is denoted by $g_E$. Whenever necessary, we will denote by
$L^2_E$ and $L^2_P$ the Sobolev norms in $\Gamma(\Lambda^*V)$ with
respect to $g_E$ and $g_P$, respectively, together with the
hermitian metric in $V$.

Model solutions of (\ref{hiteqs}) in a neighbourhood of the
singularities were described by Biquard \cite{B4}:
\begin{eqnarray*}
B & = & b\frac{d\xi}{\xi}+b^*\frac{d\overline{\xi}}{\overline{\xi}} \\
\Phi & = & \phi_0\frac{d\xi}{\xi}
\end{eqnarray*}
where $b,\phi_0\in\mathfrak{sl}(k)$. Every meromorphic Higgs pair with 
a simple pole approaches this model close enough to the singularities.

Finally, a Higgs pair $(B,\Phi)$ is said to be {\em admissible} if
$V$ has no covariantly constant sections.


\section{Construction of doubly-periodic instantons} \label{inv}

Our task now is to construct a $SU(2)$ vector bundle over $\torus$,
with an instanton connection on it, starting from a suitable singular
Higgs pair as described in the previous section.

The key feature of Nahm transforms is to try to solve a Dirac
equation, and then use its $L^2$-solutions to form a vector bundle
over the dual lattice; see the references in the introduction.

So let $S^+=\Lambda^0\oplus\Lambda^{1,1}$ and $S^-=\Lambda^{1,0}\oplus\Lambda^{0,1}$,
as vector bundles over $\dual$. The idea is to study the following
elliptic operators:
\begin{eqnarray}
{\cal D}:\Gamma(V\otimes S^+) \seta \Gamma(V\otimes S^-) & &
{\cal D}^*:\Gamma(V\otimes S^-) \seta \Gamma(V\otimes S^+) \nonumber \\
{\cal D}=(\del_B+\Phi)-(\del_B+\Phi)^* & &
{\cal D}^*=(\del_B+\Phi)^*-(\del_B+\Phi) \label{tmp.dirac}
\end{eqnarray}
where $(B,\Phi)$ is a Higgs pair. Note that the operators in
(\ref{tmp.dirac}) are just the Dirac operators coupled to the
connection $\widetilde{B}$, obtained by lifting the Higgs pair
$(B,\Phi)$ to an invariant ASD connection on $\real^4$, as
above.

The next step is to prove that the admissibility condition implies
the vanishing of the $L^2$-kernel of ${\cal D}$:

\begin{prop} \label{adm2}
The Higgs pair $(B,\Phi)$ is admissible if and only if \linebreak
$L^2_E{\rm -ker}{\cal D}=\{0\}$.
\end{prop}
\pf
Given a section $s\in L^2_2(V\otimes S^+)$, the Weitzenb\"ock
formula with respect to the Euclidean metric on the punctured torus
is given by:
\begin{eqnarray*}
(\del_B^*\del_B+\del_B \del_B^*)s & = &\nabla^*_B\nabla_Bs+F_Bs \ = \
\nabla^*_B\nabla_Bs-[\Phi,\Phi^*]s \\
\Rightarrow \ \ \ \nabla^*_B\nabla_Bs & = &
(\del_B^*\del_B+\del_B\del_B^*+\Phi\Phi^*+\Phi^*\Phi)s \\
& = & \left\{ (\del_B+\Phi)(\del_B^*+\Phi^*)+(\del_B^*+\Phi^*)(\del_B+\Phi) \right\}s \\
& = &{\cal D}^*{\cal D}s
\end{eqnarray*}
and integrating by parts, we get:
$$ ||{\cal D}s||^2_{L^2_E} \ = \ ||\nabla_Bs||^2_{L^2_E} $$
Thus, if $B$ is admissible, then the $L^2_E$-kernel of $\cal D$
must vanish. The converse statement is also clear.
\pfend

In other words, the above proposition implies that the
$L^2_E$-cohomology of orders 0 and 2 of the complex:
\begin{equation} \label{tmp.cpx}
{\cal C}\ : \
0 \seta \Lambda^0V \stackrel{\Phi+\del_B}{\longrightarrow}
  \Lambda^{1,0}V\oplus\Lambda^{0,1}V \stackrel{\del_B+\Phi}{\longrightarrow}
  \Lambda^{1,1}V \seta 0
\end{equation}
must vanish. On the other hand, since the $L^2$-norm for 1-forms is
conformally invariant, the $L^2$-cohomology $H^1({\cal C})$ does
not depend on the metric itself, only on its conformal class.

Motivated by a result of Biquard (theorem 12.1 in \cite{B3}) we
will see how one can identify $H^1({\cal C})$ in terms of a certain
hypercohomology vector space which we now introduce.

Let $\vv\seta\dual$ be the extended holomorphic vector bundle mentioned
above. Recall that if $\ksi_0$ is not an element of order 2 then the
residue of the Higgs field $\Phi$ at $\as$ is a $k\times k$ matrix of
rank 1. Therefore, if $s$ is a local holomorphic section on a neighbourhood
of $\as$, $\Phi(s)$ has at most a simple pole at $\as$ and its residue
has the form $(*,0,\ldots,0)$ on some suitable trivialisation.

Similarly, if $\ksi_0$ is an element of order 2, $\Phi(s)$ has at most
a simple pole at $\as$ and its residue has the form $(*,*,0,\ldots,0)$
on some suitable trivialisation.

This local discussion motivates the definition of a sheaf $\pp_{\as}$
such that, given an open cover $\{U_\alpha\}$ of $\dual$:
\begin{itemize}
\item $\pp_{\as}(U_\alpha)=\oo_{\dual}(\vv)(U_\alpha)$, if $\as\notin U_\alpha$;
\item $\pp_{\as}(U_\alpha)=\{$meromorphic sections of $U_\alpha\seta U_\alpha\times\cpx^k$
which have at most a simple pole at $\as$ with residue lying either along a
2-dimensional subspace of $\cpx^k$ if $\ksi_0$ has order 2, or along a 1-dimensional
subspace of $\cpx^k$ otherwise$\}$, if $\as\in U_\alpha$.
\end{itemize}
It is easy to see that such $\pp_{\as}$ is a coherent sheaf. To simplify notation,
we drop the subscript $\as$ out.

Hence, $\Phi$ can be regarded as the map of sheaves:
\begin{equation} \label{hypcpx}
\Phi : \vv \seta \pp\otimes K_{\dual}
\end{equation}
Seen as a two-term complex of sheaves, the map (\ref{hypcpx}) induces
an exact sequences of hypercohomology vector spaces:
\begin{equation} \label{hcoho} \begin{array}{cccccccc}
0 & \seta & {\Bbb H}^0(\dual,\Phi) & \seta & H^0(\dual,\vv) &
    \stackrel{\Phi}{\seta} & H^0(\dual,\pp\otimes K_{\dual}) & \seta \\
  & \seta & {\Bbb H}^1(\dual,\Phi) & \seta & H^1(\dual,\vv) &
    \stackrel{\Phi}{\seta} & H^1(\dual,\pp\otimes K_{\dual}) & \seta \\
  & \seta & {\Bbb H}^2(\dual,\Phi) & \seta & 0 & & &
\end{array} \end{equation}

It is easy to see that:
\begin{eqnarray*}
{\Bbb H}^0(\dual,\Phi) & = & {\rm ker}\left\{ H^0(\dual,\vv)
\stackrel{\Phi}{\seta} H^0(\dual,\pp\otimes K_{\dual}) \right\} \\
{\Bbb H}^2(\dual,\Phi) & = & {\rm coker}\left\{ H^1(\dual,\vv)
\stackrel{\Phi}{\seta} H^1(\dual,\pp\otimes K_{\dual}) \right\}
\end{eqnarray*}
and admissibility implies that the right-hand sides must vanish:
restricted to $\dual\setminus\{\as\}$, a section there would give a
section in the kernel of ${\cal D}$ (or, equivalently, a class in
$H^0({\cal C})$ and $H^1({\cal C})$). Therefore, the dimension of
${\Bbb H}^1(\dual,\Phi)$ is equal to $\chi(\pp\otimes
K_{\dual})-\chi(\vv)=\chi(\pp)-\chi(\vv)$.

To compute this number, note that there is also a natural map \linebreak
$\vv\stackrel{\iota}{\seta}\pp$ defined as the {\em local inclusion}
of holomorphic local sections (elements of $\oo_{\dual}(\vv)(U_\alpha)$),
into the meromorphic ones (elements of $\pp(U_\alpha)$). It fits
into the following sequence of sheaves:
\begin{eqnarray}
0\seta \vv \stackrel{\iota}{\seta} \pp
\stackrel{res_{\ksi_0}}{\longrightarrow} \rr_{\ksi_0}\seta 0
& & {\rm if}\ \ksi_0\ {\rm has\ order\ 2,} \label{id1} \\
0\seta \vv \stackrel{\iota}{\seta} \pp
\stackrel{res_{\as}}{\longrightarrow} \rr_{\as}\seta 0 & &
{\rm otherwise} \label{id2}
\end{eqnarray}
where $\rr_{\ksi_0}$ is the skyscraper sheaf supported at $\ksi_0$
and stalk isomorphic to $\cpx^2$ and $\rr_{\as}$ is the skyscraper
sheaf supported at $\as$ and stalks isomorphic to $\cpx$. Since
$\chi(\rr_{\as})=\chi(\rr_{\xi_0})=2$, we conclude that
${\Bbb H}^1(\dual,\Phi)$ is a 2-dimensional complex vector space.

\begin{prop} \label{ker/hcoh}
The hypercohomology induced by the map of sheaves (\ref{hypcpx})
coincides with the $L^2_P$-cohomology of the complex (\ref{tmp.cpx}).
\end{prop}

In particular, we have identifications:
$$ {\Bbb H}^1(\dual,\Phi) \ \equiv \ L^2_P{\rm -cohomology}\ H^1({\cal C})
   \ \equiv \ L^2_E{\rm -cohomology}\ H^1({\cal C}) $$
Furthermore, note also that the $L^2_E$-cohomology of 1-forms with
respect to the Euclidean metric is a 2-dimensional complex vector
space.

\pf The hypercohomology defined by the map (\ref{hypcpx}) is given by the
total cohomology of the double complex:
\begin{eqnarray*}
\Lambda^0\vv & \stackrel{\Phi}{\seta} & \Lambda^{1,0}\pp \\
\del\ \downarrow & & \downarrow\ \del \\
\Lambda^{0,1}\vv & \stackrel{\Phi}{\seta} & \Lambda^{1,0}\pp
\end{eqnarray*}
which in turns is just the cohomology of the complex:
$$ 0\seta \Lambda^0\vv \stackrel{\Phi+\del}{\seta} \Lambda^{1,0}\pp \oplus
\Lambda^{0,1}\vv \stackrel{\del+\Phi}{\seta} \Lambda^{1,0}\pp \seta 0 $$
Now restricting the complex above to the punctured torus
$\dual\setminus\{\as\}$, we get:
$$ 0 \seta \Lambda^0V \stackrel{\Phi+\del_B}{\seta} \Lambda^1V
\stackrel{\del_B+\Phi}{\seta} \Lambda^2V \seta0 $$
which is, of course, the complex ${\cal C}$.

So, let $s$ be a section of $\Lambda^{1,0}\pp\oplus\Lambda^{0,1}\vv$
defining a class in ${\Bbb H}^1(\dual,\Phi)$. Thus, restricting $s$
to $\dual\setminus\{\as\}$ yields a section $s_r$ of $L^2(\Lambda^1V)$
defining a class in $H^1({\cal C})$.

Such {\em restriction map} is clearly a well-defined map:
\begin{eqnarray*}
R:{\Bbb H}^1(\dual,\Phi) & \seta & H^1({\cal C}) \\
<s> & \seta & <s_r>
\end{eqnarray*}
We claim that it is also injective. Indeed, suppose that $s_r$
represents the zero class, i.e. there is $t\in L^2_2(\Lambda^0V)$
such that $s_r=(\del_B+\Phi)t$. However, $L^2_2\hookrightarrow C^0$
is a bounded inclusion in real dimension 2. Thus, $h(t,t)$ must be
bounded at the punctures $\as$, and $t$ must be itself bounded
along the kernel of the residues of $\Phi$. On the other hand, the
hermitian metric degenerates along the image of the residues of
$\Phi$, so $t$ might be singular on this direction. Indeed, 
$h\sim O(r^{1\pm\alpha})$ in a holomorphic trivialisation, so that 
$t\sim O(r^{-\frac{1}{2}(1\pm\alpha)})$. But then the derivatives 
of $t$ will not be square integrable, contradicting our hypothesis 
that $t$ belongs to $L^2_2$. So $t$ must be bounded at $\as$.

This implies that $t\in L^2_2(\Lambda^0\vv)$ also with respect to
the $h'$ metric, so that $s_r$ is indeed the restriction of a
section representing the zero class in ${\Bbb H}^1(\dual,\Phi)$.

Finally, to show that $R$ is an isomorphism, it is enough by admissibility
to argue that the $L^2$ index of the complex ${\cal C}$ is $-2$.

It was shown by Biquard (theorem 5.1 in \cite{B3}) the laplacian
associated to the complex ${\cal C}$ is Fredholm when acting
between $L^2_P$ sections. This implies that ${\cal D}$ is also
Fredholm. Its index can be computed via Gromov-Lawson's relative
index theorem, and it coincides with the index of the Dirac
operator on $\vv$:
$$ \ind({\cal D})=\ind(\del_B-\del_B^*)={\rm deg}\vv=-2 $$
as desired
\pfend

\paragraph{Constructing the transformed bundle.}
We are finally in a position to construct a vector bundle with
connection over $\torus$ out of a Higgs pair $(B,\Phi)$. Recall that
$\jj(\dual)=T$, the {\em Jacobian} of $\dual$, is defined as the set
of flat holomorphic line bundles over $\dual$. Each $z\in T$ corresponds
to a flat holomorphic line bundle $L_z\seta\dual$. Moreover, $T$ and
$\dual$ are isomorphic as elliptic curves.

These line bundles can be given a natural constant connection compatible
with the holomorphic structure. This follows from the differential-geometric
definition of $T$:
$$ T = \{z\in(\real^2)^*\ |\ z(\xi)\in\zed,\forall \xi\in\Lambda \} $$
where $\Lambda\subset\real^4$ is the two-dimensional lattice
generating $\dual$. Hence each $z\in T$ can be regarded as a
constant, real 1-form over $\dual$, so that $\omega_z=i\cdot z$ is a
connection on a topologically trivial line bundle $L\seta\dual$. Each
such connection defines a different holomorphic structure on $L$,
which we denote by $L_z$.

Conversely, $\dual$ parametrises the set of holomorphic flat line
bundles with connection over $T$. Each point $\xi\in\dual$ corresponds
to the line bundle $L_\ksi\seta T$ with a connection $\omega_\ksi$.

Now consider the restrictions $L_z\seta\dual\setminus\{\as\}$, with
its natural connection $\omega_z$, and form the tensor product
$V(z)=V\otimes L_z$. The connection $B$ can be tensored with
$\omega_z$ to obtain another connection that we denote by $B_z$.

Let $i:V(z)\seta V(z)$ be the identity bundle automorphism and define
$\Phi_{w}=\Phi-w\cdot i$, where $w$ is a complex number. It is easy
to see that $(B_z,\Phi_w)$ is still an admissible Higgs pair, for all
$(z,w)\in\torus$.

Therefore, we get the following continuous family of Dirac-type operators:
\begin{equation} \label{prt}
{\cal D}_{(z,w)}=(\del_{B_z}+\Phi_w)-(\del_{B_z}+\Phi_w)^*
\end{equation}
From proposition \ref{adm2}, we have that $L^2_E{\rm -ker}{\cal D}_{(z,w)}$
vanishes for all $(z,w)\in\torus$. Since its index remains invariant
under this continuous deformation, we conclude that
$L^2_E{\rm -ker}{\cal D}_{(z,w)}^*$  has constant dimension equal to 2.

Define a trivial Hilbert bundle $H\seta\torus$ with fibres given by
\linebreak $L^2(V(z)\otimes S^-)$. It follows that 
$E_{(z,w)}={\rm ker}{\cal D}_{(z,w)}^*$ forms a vector sub-bundle
$E\stackrel{i}{\hookrightarrow}H$ of rank 2. Furthermore \cite{DK},
$E$ is also equipped with an hermitian metric, induced from the
$L^2$ metric on $H$, and an unitary connection $A$, defined as
follows:
\begin{equation} \label{invconn}
\nabla_A=P \circ d \circ i
\end{equation}
where $d$ means differentiation with respect to $(z,w)$ on the
trivial Hilbert bundle (i.e. the trivial product connection) and
$P$ is the fibrewise orthogonal projection
$P:L^2(V(z)\otimes S^-)\seta {\rm ker}{\cal D}_{(z,w)}^*$.
Clearly, $A$ defined on (\ref{invconn}) is unitary.

Note also that the hermitian metric in $H$ is actually
conformally invariant with respect to the choice of metric in
$\dual\setminus\{\as\}$, since the inner product in
$L^2(V(z)\otimes S^-)$ is. Therefore, the induced hermitian
metric in $E$ is also conformally invariant.

\paragraph{Monad description.}
The transformed bundle $E$ also admits a monad-type description. More
precisely, once a metric is chosen, the family of Dirac operators
${\rm ker}{\cal D}_{(z,w)}^*$ can be unfolded into the following
family of elliptic complexes ${\cal C}(z,w)$:
\small \baselineskip18pt
\begin{equation} \label{mnd}
0 \seta L^2_{2,E}(\Lambda^0V(z)) \stackrel{\Phi_w+\del_{B_z}}{\longrightarrow}
L^2_{1,E}(\Lambda^{1,0}V(z)\oplus\Lambda^{0,1}V(z)) \stackrel{\del_{B_z}+\Phi_w}{\longrightarrow} L^2_E(\Lambda^{1,1}V(z)) \seta 0
\end{equation}
\normalsize \baselineskip18pt

Admissibility implies that $H^0({\cal C}(z,w))$ and $H^2({\cal C}(z,w))$
must vanish, and $H^1({\cal C}(z,w))$ coincides with
$L^2_E{\rm -ker}{\cal D}^*_{(z,w)}$. As $(z,w)$ sweeps out $\torus$,
$H^1({\cal C}(z,w))$ forms a rank 2 holomorphic vector bundle with a
natural hermitian metric and a compatible unitary connection $A$,
equivalent to the ones defined as above; see \cite{DK}.

\subsection{Anti-self-duality and curvature decay}
The next proposition fulfills the first goal of this paper,
i.e. to show that the connection $A$ defined above is in fact a
finite-energy anti-self-dual instanton on the rank 2 bundle
$E\seta\torus$. We say $f\sim O(|w|^n)$ if the complex function
$f:\cpx\seta\cpx$ satisfies:
\begin{equation} \label{defn.o}
\lim_{|w|\seta\infty}\frac{|f(w)|}{|w|^n}<\infty
\end{equation}

\begin{thm} \label{asd.inv}
The transformed connection $A$ is anti-self-dual with respect
to the Euclidean metric. Furthermore, its curvature satisfies
$|F_A|\sim O(|w|^{-2})$.
\end{thm}

\pf
Since $A$ is an unitary connection, we only have to verify
that the component of $F_A$ along the K\"ahler class $\kappa$
of $\torus$ vanishes.

Let $\{\psi_1,\psi_2\}$ be a local holomorphic frame for $E$,
orthonormal with respect to the hermitian metric induced from $H$.
Fix some $(z,w)\in\torus$ so that, as a section of $\vv(z)\otimes S^-\seta\dual$,
we have $\psi_i=\psi_i(\ksi;z,w)\in{\rm ker}{\cal D}_{(z,w)}^*$.

In this trivialisation, the matrix elements of the
curvature $F_A$ can then be written as follows:
\begin{eqnarray}
(F_A)_{ij} & = & \langle \psi_j,\nabla_A\nabla_A\psi_i \rangle
\ = \ \langle \psi_j,d\circ P\circ d\psi_i \rangle \ = \nonumber \\
& = & \langle {\cal D}_{(z,w)}^*(d\psi_j),G_{(z,w)}{\cal D}_{(z,w)}^*(d\psi_j) \rangle
\end{eqnarray}
where the inner product is taken in $L^2(V(z)\otimes S^-)$, integrating
out the $\ksi$ coordinate; the finiteness of the integral is guaranteed 
by the fact that $\psi_j\in L^2_1(V(z)\otimes S^-)$. Note also that the
inner product is conformally invariant with respect to the choice of 
metric on $\dual\setminus\{\as\}$. Hence, the expression for the curvature
above is the same for both the Euclidean and Poincar\'e metrics.

Moreover, $G_{(z,w)}$ is the Green's operator for
${\cal D}_{(z,w)}^*{\cal D}_{(z,w)}$. Note that:
$$ [{\cal D}_{(z,w)}^*,d]\psi_i=\Omega^\prime\cdot \psi_i $$
where $\Omega^\prime=(idz_1+dw_1)\wedge d\ksi_1+(idz_2+dw_2)\wedge
d\ksi_2$ and ``$\cdot$'' denotes Clifford multiplication. So,
\begin{eqnarray}
\kappa\llcorner(F_A)_{ij} & = & \langle \psi_j,
\underbrace{\kappa\llcorner(\Omega'\wedge\Omega')}\cdot G_{(z,w)}\psi_i \rangle
\ = \ 0 \label{curv} \\
& & \hspace{1.5cm} =0 \nonumber
\end{eqnarray}
and this proves the first statement.

It is easy to see from (\ref{curv}) that the asymptotic behaviour of
$|(F_A)_{ij}|$ depends only on the behaviour of the operator norm
$||G_{(z,w)}||$ for large $|w|$.

We can estimate $||G_{(z,w)}||$ by looking for a lower bound for the
eigenvalues of the associated laplacian acting on $V\otimes S^-$:
\begin{equation} \label{zwlap}
{\cal D}_{(z,w)}{\cal D}_{(z,w)}^* \ = \ {\cal D}_{z}{\cal D}_{z}^* -
w\phi^* - \overline{w}\phi+|w|^2
\end{equation}
where ${\cal D}_{z}={\cal D}_{(z,w=0)}$ and $\Phi=\phi d\ksi$, with
$\phi\in{\rm End}V$; $\phi^*$ denotes the adjoint (conjugate transpose)
endomorphism.

In other words, we want to find a lower bound for the following expression:
\begin{eqnarray}
& \left| \langle ({\cal D}_{z}{\cal D}_{z}^* + |w|^2)s,s \rangle -
         \langle(w\phi^* + \overline{w}\phi)s,s \rangle \right| \geq & \nonumber \\
& \geq \left| \  \langle({\cal D}_{z}{\cal D}_{z}^* + |w|^2)s,s  \rangle \ - \
            | \langle (w\phi^* + \overline{w}\phi)s,s \rangle | \ \right| \label{fin} &
\end{eqnarray}
for $s\in L^2(V\otimes S^-)$ of unit norm.

For the first term in the second line, it is easy to see that:
\begin{equation} \label{w2term}
|  \langle ({\cal D}_{z}{\cal D}_{z}^*+|w|^2)s,s  \rangle | \ = \
||{\cal D}_{z}^*s||^2 + |w|^2\cdot||s||^2 \ = \ c_1 + |w|^2
\end{equation}
for some non-zero constant $c_1=||{\cal D}_{z}^*||^2$ depending only
on $z\in T$.

The second term in (\ref{fin}) is more problematic; first note
that:
$$ |  \langle (w\phi^* + \overline{w}\phi)s,s  \rangle | \ \leq \
   |w| \cdot \left( | \langle \phi(s),s \rangle | + | \langle \phi^*(s),s \rangle | \right) $$
In a small neighbourhood $D_0$ of each singularity $\as$, we have:
\begin{eqnarray*}
\langle \phi(s),s \rangle_{L^2(D_0)} & = &
\int_{D_0} \langle \frac{\phi_0(s)}{\ksi},s \rangle rdrd\theta+
\left( \begin{array}{c} {\rm regular} \\ {\rm terms} \end{array} \right) \\
& \sim & \int_{D_0} \frac{|\phi_0|}{r}\cdot|s|^2 rdrd\theta
+ \left(\begin{array}{c} {\rm regular} \\ {\rm terms} \end{array}\right)
\end{eqnarray*}
Let $1<p<2$; using H\"older inequality, we obtain:
\begin{eqnarray*}
\int_{D_0} \frac{|\phi_0|}{\ksi}\cdot|s|^2 & \leq &
\left\{ \int_{D_0} \left( \frac{|\phi_0|}{r} \right)^p rdrd\theta \right\}^{1/p}
\left\{ \int_{D_0} |s|^{2q} \right\}^{1/q} \\
& \leq & c \cdot ||s||^2_{L^{2q}}
\end{eqnarray*}
where $q=\frac{p}{p-1}$, and for some real constant $c$ depending
only on $\phi_0$ and on the choice of $p$.

Since $2q>4$, the Sobolev embedding theorem tells us that
$L^2_1\hookrightarrow L^{2q}$ is a bounded inclusion (in real
dimension 2). In other words, there is a constant $C$ depending
only on $q$ such that $||s||_{L^{2q}} \leq C \cdot ||s||_{L^2_1}$.
Thus, arguing similarly for the $\langle \phi^*(s),s \rangle$ term,
we conclude that:
$$ |  \langle (w\phi^* + \overline{w}\phi)s,s  \rangle | \ \leq \ c_2 \cdot |w| $$
where $c_2$ is a real constant depending neither on $z$ nor on $w$, but only on
the Higgs field itself and on the choice of $p$.

Putting everything together, we have:
$$ \left| \langle ({\cal D}_{z}{\cal D}_{z}^* - w\phi^* - \overline{w}\phi + |w|^2)s,
   s \rangle \right| \geq \left| |w|^2-c_2|w|+c_1 \right| $$
so that
$$ \lim_{|w|\seta\infty}|w|^2\cdot||G_{(z,w)}||<1 $$
and the statement follows.
\pfend

\noindent {\bf Remark 1:} Note in particular that
$F_A\in L^2(\Lambda^2\otimes E)$ with respect to the Euclidean metric
on $\torus$, coming from the quotient $(\real^4)^*/\Lambda^*$. This
concludes our first task.

\noindent {\bf Remark 2:} It is also not difficult to see that gauge
equivalent Higgs pairs $(B,\Phi)$ and $(B',\Phi')$ will produce
gauge equivalent instantons $A$ and $A'$. The dependence of $A$ on
the Higgs pair $(B,\Phi)$ is contained on the $L^2$-projection operator
$P$, that is on the two linearly independent solutions of \linebreak
${\cal D}^*_{(z,w)}\psi=0$. Gauge equivalence of $(B,\Phi)$ and
$(B',\Phi')$ gives an automorphism of the transformed bundle $E$,
in other words, a gauge equivalence between $A$ and $A'$.

\noindent {\bf Remark 3:} The instanton connection $A$ induces a
holomorphic structure $\del_A$ on the the transformed bundle
$E\seta\torus$.

In order to further understand the asymptotic behaviour of the transformed
connection,  we must now pass to an equivalent holomorphic description
of the above transform.


\section{Holomorphic version and extensibility}
\label{invholo}

Motivated by curvature decay established above, one can expect to
find a holomorphic vector bundle $\ee\seta\tproj$ which extends
$(E,\del_A)$. The idea is to find a suitable perturbation of the
Higgs field $\Phi$ for which $w=\infty$ makes sense.

As above, the torus parameter $z\in T$ simply twists the holomorphic
bundle $\vv\seta\dual$. We denote:
\begin{equation} \begin{array}{ccc}
\vv(z)=\vv\otimes L_z & \ \ \ & \pp(z)=\pp\otimes L_z
\end{array} \end{equation}
Since $\Phi\in H^0(\dual,{\rm Hom}(\vv,\pp)\otimes K_{\dual})$,
tensoring both sides of (\ref{hypcpx}) by the line bundle $L_z$
does not alter the sheaf homomorphism $\Phi$, so we have the
family of maps:
$$ \Phi : \vv(z) \seta \pp(z)\otimes K_{\dual} $$
parametrised by $z\in T$.

To define the perturbation $\Phi_w$, recall that, regarding
$\proj=\cpx\cup\{\infty\}$, we can fix two holomorphic sections
$s_0,s_\infty\in H^0(\proj,\oo_{\proj}(1))$ such that $s_0$ vanishes
at $0\in\cpx$ and $s_\infty$ vanishes at the point added at infinity.
In homogeneous coordinates $\{(w_1,w_2)\in\cpx^2|w_2\neq0\}$ and
$\{(w_1,w_2)\in\cpx^2|w_1\neq0\}$, we have that, respectively
($w=w_1/w_2$):
\begin{eqnarray*}
s_0(w)=w & \ \ \ & s_0(w)=1 \\
s_\infty(w)=1 & \ \ \ & s_\infty(w)=\frac{1}{w}
\end{eqnarray*}

Consider now the map of sheaves parametrised by pairs $(z,w)\in\tproj$:
\begin{eqnarray}
& \Phi_w: \vv(z) \seta \pp(z)\otimes K_{\dual} & \nonumber \\
& \Phi_w = s_\infty(w)\cdot\Phi-s_0(w)\cdot\iota\cdot d\ksi & \label{defhiggs}
\end{eqnarray}
Clearly, on $\proj\setminus\{\infty\}=\cpx$ this is just
$\Phi_w=\Phi-w\cdot\iota$, the same perturbation we defined before.
Moreover, if $w=\infty$, then $\Phi_\infty = \iota\cdot d\ksi$

The hypercohomology vector spaces ${\Bbb H}^0(\dual,\Phi_w)$ and
${\Bbb H}^2(\dual,\Phi_w)$ of the two-term complex (\ref{defhiggs})
must vanish by admissibility. On the other hand,
${\Bbb H}^1(\dual,\Phi_w)$ also makes sense for $\infty\in\proj$,
and we can define a $SU(2)$ holomorphic vector bundle
$\ee\seta\tproj$ with fibres given by
$\ee_{(z,w)}={\Bbb H}^1(\dual,\Phi_w)$. Moreover, $\ee$ is
actually a {\em holomorphic extension} of $(E,\del_A)$, in the
sense that, holomorphically:
\begin{equation} \label{holoextn}
\ee|_{T\times(\proj\setminus\{\infty\})} \simeq (E,\del_A)
\end{equation}

Equivalently, $\ee$ can be seen as the hermitian holomorphic vector
bundle induced by the monad
\begin{equation} \label{op.db.cpx}
0\seta \Lambda^0\vv \stackrel{\Phi+\del}{\seta} \Lambda^{1,0}\pp \oplus
\Lambda^{0,1}\vv \stackrel{\del+\Phi}{\seta} \Lambda^{1,0}\pp \seta 0
\end{equation}

Consider the metric $H'$ induced from the monad (\ref{op.db.cpx}) above,
while $H$ is induced from the monad (\ref{mnd}). Now, $H$ is bounded above
by $H'$ because the hermitian  metric $h$ on the bundle $V$ in (\ref{mnd})
is bounded above by the metric $h'$ on the bundle $\vv$ in (\ref{op.db.cpx}).

We now show that the position of the singularities of the Higgs
pair determines the holomorphic type of the restriction of the
extended transformed bundle over the added divisor at infinity.
First, recall that there is an unique line bundle ${\bf P}\seta
T\times\dual$, the so-called {\em Poincar\'e line bundle},
satisfying:
\begin{eqnarray*}
{\mathbf P}|_{T\times\{\ksi\}}\simeq L_\ksi & \ \ &
{\mathbf P}|_{\{z\}\times\dual}\simeq L_{-z}
\end{eqnarray*}
It can be constructed as follows. Identifying $T$ and $\dual$ as before,
let $\Delta$ be the diagonal inside $T\times\dual$, and consider the
divisor $D=\Delta-T\times{\hat{e}}-{e}\times\dual$. Then
${\mathbf P}=\oo_{T\times\dual}(D)$; it is easy to see that the sheaf so
defined restricts as wanted.

Note that although the two restrictions above are flat line bundles over $T$ and
$\dual$ respectively, the Poincar\'e bundle itself is not topologically trivial;
in fact, $c_1({\mathbf P})\in H^1(T)\otimes H^1(\dual)\subset H^2(T\times\dual)$.
More precisely, the unitary connection and its corresponding curvature are given by:
$$ \omega(z,\ksi)=i\pi\cdot\sum_{\mu=1}^{2} \big(\ksi_\mu dz_\mu - z_\mu d\xi_\mu\big)
   \ \ {\rm and} \ \
   \Omega(z,\ksi)=2i\pi\cdot\sum_{\mu=1}^{2}d\ksi_\mu\wedge dz_\mu $$
Restricting to each $T\times\{\ksi\}$, the line bundles $L_\ksi\seta T$ are given
flat connections \linebreak $\omega_\ksi=i\pi\cdot\sum_{\mu=1}^{2}\ksi_\mu dz_\mu$,
with constant coefficients. Similarly, the line bundles $L_z\seta\dual$ are given 
the flat connections $\omega_z=-i\pi\cdot\sum_{\mu=1}^{2}z_\mu d\xi_\mu$ as
described in the previous section.
Finally, note that:
$$ c_1({\bf P}) = \frac{i}{2\pi}\Omega \ \ \Rightarrow \ \ 
   c_1({\bf P})^2 = -2\cdot t\wedge\hat{t} $$ 
where $t$ and $\hat{t}$ are the generators of $H^2(T)$ and $H^2(\dual)$,
respectively.

\begin{lem} \label{infty}
$\ee|_{T_\infty}\equiv L_{\ksi_0}\oplus L_{-\ksi_0}$
\end{lem}

\pf Substituting $w=\infty\in\proj$, we get from (\ref{defhiggs})
that $\Phi_\infty=\iota\cdot d\ksi$. Therefore, the induced
hypercohomology sequence (\ref{hcoho2}) coincides with the long
exact sequence of cohomology induced by the sheaf sequences
(\ref{id1}) and (\ref{id2}),
which is given by:
\begin{equation} \label{hcoho4} \begin{array}{rcl}
0 & \seta & H^0(\dual,\vv(z)) \stackrel{\Phi_\infty}{\seta}
            H^0(\dual,\pp(z)\otimes K_{\dual}) \seta H^0(\dual,\rr_{\as}(z)) \seta \\
  & \seta & H^1(\dual,\vv(z)) \stackrel{\Phi_\infty}{\seta}
            H^1(\dual,\pp(z)\otimes K_{\dual}) \seta 0
\end{array} \end{equation}
Hence, ${\Bbb H}^1(\dual,(z,\infty))=H^0(\dual,\rr_{\as}(z))$.
The right hand side is canonically identified with
$(L_z)_{\ksi_0}\oplus (L_z)_{-\ksi_0}$, where by $(L_z)_{\ksi_0}$
we mean the fibre of $L_z\seta\dual$ over the point $\ksi_0\in\dual$.

On the other hand, $(L_z)_{\ksi_0}={\bf P}_{(z,\ksi_0)}=(L_{\ksi_0})_z$,
where ${\bf P}\seta T\times\dual$ is the Poincar\'e line bundle. Thus,
the bundle over $T_\infty$ with fibres given by $H^0(\dual,\rr_{\as}(z))$
is isomorphic to $L_{\ksi_0}\oplus L_{-\ksi_0}$, as we wished to prove.
\pfend

The topological type of $\ee$ is also fixed from the initial data:
the rank of the bundle $V$ is translated into the second Chern class
of the extended transformed bundle $\ee$. In the next lemma, we denote
the generator of $H^2(\proj,\zed)$ by $p$.

\begin{lem} \label{chch}
$ch(\ee) = 2 - k\cdot t\wedge p$.
\end{lem}
\pf
The exact sequence:
\begin{equation} \label{hcoho2} \begin{array}{ccl}
0 & \seta & H^0(\dual,\vv(z)) \stackrel{\Phi_w}{\seta}
            H^0(\dual,\pp(z)\otimes K_{\dual}) \seta {\Bbb H}^1(\dual,(z,w)) \seta \\
  & \seta & H^1(\dual,\vv(z)) \stackrel{\Phi_w}{\seta}
            H^1(\dual,\pp(z)\otimes K_{\dual}) \seta 0
\end{array} \end{equation}
induces a sequence of coherent sheaves over $\torus$, with stalks over
$(z,w)$ given by the above cohomology groups:
\begin{equation} \label{hcoho3} \begin{array}{ccl}
0 & \seta & {\cal H}^0(\dual,\vv(z)) \stackrel{\Phi_w}{\seta}
            {\cal H}^0(\dual,\pp(z)\otimes K_{\dual}) \seta \check{\ee} \seta \\
  & \seta & {\cal H}^1(\dual,\vv(z)) \stackrel{\Phi_w}{\seta}
            {\cal H}^1(\dual,\pp(z)\otimes K_{\dual}) \seta 0
\end{array} \end{equation}
In this way, the Chern character of $\check{\ee}$ will then be given
by the alternating sum of the Chern characters of these sheaves, which can
be computed via the usual Grothendieck-Riemann-Roch for families.

Consider the bundle ${\bf G}_1\seta\tproj\times\dual$ given by
${\bf G}_1=p_3^*\vv\otimes p_{13}^*{\bf P}$. Clearly,
${\bf G}_1|_{(z,w)\times\dual}=\vv(z)$, so that:
\begin{equation} \label{yyy}
ch({\cal H}^0(\dual,\vv(z)))-ch({\cal H}^1(\dual,\vv(z)))=
ch({\bf G}_1)td(\dual)/[\dual]
\end{equation}

Now consider the sheaf:
${\bf G}_2=p_3^*{\cal P}\otimes p_{13}^*{\bf P}\otimes p_2^*\oo_{\proj}(1)$.
The twisting by $\oo_{\proj}(1)$ accounts for the multiplication by the
section $s_0\in H^0(\proj,\oo_{\proj}(1))$ contained in $\Phi_w$. As
above, ${\bf G}_1|_{(z,w)\times\dual}=\pp(z)$, and we have:
\begin{equation} \label{yyz}
ch({\cal H}^0(\dual,\pp(z)\otimes K_{\dual})) -
ch({\cal H}^1(\dual,\pp(z)\otimes K_{\dual})) = ch({\bf G}_2)td(\dual)/[\dual]
\end{equation}
Therefore:
\begin{eqnarray*}
ch(\ee) & = & (\ref{yyz})-(\ref{yyy}) \ = \\ & = & 
\left( c_1(\pp) - c_1(\vv) + c_1(\pp)\wedge p + \frac{k}{2}c_1({\bf P})^2\wedge p \right)/
[\dual] \ = \\
& = & \chi(\pp) - {\rm deg}\vv + \chi(\pp)\cdot p - k\cdot t\wedge p \ = \
2 - k\cdot t\wedge p
\end{eqnarray*}
as desired.
\pfend

Finally, we argue that the determinant bundle of $\ee$ is trivial,
so that $A$ is indeed a $SU(2)$ instanton. Note that ${\rm det}\ee$
is a line bundle with vanishing first Chern class, so it must be
the pull back of a flat line bundle $L_\ksi\seta T$. But
${\rm det}\ee|_{T_\infty}=\underline{\cpx}$, hence ${\rm det}\ee$
must be holomorphically trivial, as desired.

We call $\ksi_0\in\jj(T)$ the {\em asymptotic state} associated to
the doubly-periodic instanton connection $A$, and the integer $k$
its {\em instanton number}. The Nahm transform constructed above
guarantees the existence of doubly-periodic instantons of any given
charge and asymptotic state.

\subsection{Extensible doubly-periodic instanton connections}

Motivated by the properties established above, we say that a doubly-periodic
instanton connection $A$ on a bundle $E\seta\torus$ is {\em extensible} if
the following hypothesis hold:
\begin{enumerate}
\item $|F_A|\sim O(|w|^{-2})$;
\item there is a holomorphic vector bundle $\ee\seta\tproj$ with trivial
determinant such that $\ee|_{T\times(\proj\setminus\{\infty\})}\simeq(E,\del_A)$,
where $\del_A$ is the holomorphic structure on $E$ induced by the instanton
connection $A$;
\end{enumerate}

This definition will be our starting point in \cite{J2}, where we shall present
the Nahm transform of doubly-periodic instantons, i.e. the inverse of the
construction shown here.


\section{Conclusion}

In this paper we have shown how finite energy, doubly-periodic
instantons can be produced by performing a Nahm transform on
certain singular Higgs pairs. The rank of the Higgs bundle is
translated into the instanton number; the number of singularities
of the Higgs field (i.e. the degree of the holomorphic Higgs bundle
$\vv$) gives the rank of the transformed instanton, and its
positions determine how the instanton connection ``splits at
infinity''. Indeed, it is easy to generalise the above construction
by allowing more than two singularities on the original Higgs
field, so that higher rank doubly-periodic instantons are obtained;
see \cite{J2}.

Moreover, one would also like to understand how the parabolic
parameters $(\alpha,\epsilon)$ are translated into the
doubly-periodic instantons produced via the Nahm transform as
above. On general grounds, we expect these parameters to be
translated into more detailed information on the asymptotic
behaviour of $A$.

From the more analytical point of view, it is also interesting to
ask if the curvature decay (proposition \ref{asd.inv}) is enough to
ensure extensibility. More precisely, one can expect to be able to
prove the following result:

\begin{conj} \label{dec.extn}
If $A$ is anti-self-dual and $|F_A|\sim O(|w|^{-2})$, then there
is a holomorphic vector bundle $\ee\seta\tproj$ such that
$$ \ee|_{T\times(\proj\setminus\{\infty\})}\simeq(E,\del_A) $$
In other words, $A$ is extensible.
\end{conj}

Such conjecture motivates other questions:
\begin{itemize}
\item Do all anti-self-dual connections on $E\seta\torus$ with
      finite energy with respect to the Euclidean metric satisfy
      $|F_A|\sim O(|w|^{-2})$?
\item Does the converse holds, i.e. if $A$ is extensible then
      $|F_A|\sim O(|w|^{-2})$? If not, what are the necessary and
      sufficient analytical conditions for extensibility (in terms of
      the Euclidean metric)?
\item Given a holomorphic bundle $\ee\seta\tproj$, is there a
      connection $A$ on $\ee|_{T\times(\proj\setminus\{\infty\})}$
      such that $A$ is anti-self-dual and $|F_A|\sim O(|w|^{-2})$
      with respect to the Euclidean metric?
\end{itemize}
We hope to address these issues in a future paper \cite{BJ}.


 \end{document}